\documentclass[12pt]{article}
\usepackage{amsmath,amssymb,amsfonts,amsthm,url,xspace,graphicx,psfrag}
\usepackage{picinpar}
\def\volno{18}
\def\volyear{2011}
\def\papno{P130}
\usepackage{e-jc}
\usepackage[pdftitle={Double-dimer pairings and skew Young diagrams},
            pdfauthor={Richard W. Kenyon and David B. Wilson},
            bookmarks=true,bookmarksopen=true,bookmarksopenlevel=2]{hyperref}
\usepackage[all]{hypcap}

\newcommand{\sref}[1]{\hyperref[#1]{\S~\ref*{#1}}}
\newcommand{\aref}[1]{\hyperref[#1]{Appendix~\ref*{#1}}}
\newcommand{\lref}[1]{\hyperref[#1]{Lemma~\ref*{#1}}}
\newcommand{\tref}[1]{\hyperref[#1]{Theorem~\ref*{#1}}}
\newcommand{\cref}[1]{\hyperref[#1]{Corollary~\ref*{#1}}}
\newcommand{\fref}[1]{\hyperref[#1]{Figure~\ref*{#1}}}
\newcommand{\pref}[1]{\hyperref[#1]{Proposition~\ref*{#1}}}

\def\clap#1{\hbox to 0pt{\hss#1\hss}}

 \newcommand{\MRhref}[2]{\href{http://www.ams.org/mathscinet-getitem?mr=#1}{MR#2}}

\makeatletter
\def\@strippedMR{}
\def\@scanforMR#1#2#3\endscan{%
  \ifx#1M\ifx#2R\def\@strippedMR{#3}%
  \else\def\@strippedMR{#1#2#3}%
  \fi\fi}
\def\@rst #1 #2other{#1}
\newcommand\MR[1]{\relax\ifhmode\unskip\spacefactor3000 \space\fi
  \@scanforMR#1\endscan
  \MRhref{\expandafter\@rst \@strippedMR other}{\@strippedMR}}
\newcommand\MRs[1]{\relax\ifhmode\unskip\spacefactor3000 \space\fi
  \@scanforMR#1\endscan
  \MRhref{\@strippedMR}{\@strippedMR}}
\makeatother

\newtheorem{theorem}{Theorem}
\numberwithin{theorem}{section}
\newtheorem{conjecture}{Conjecture}
\newtheorem{lemma}[theorem]{Lemma}
\newtheorem{proposition}[theorem]{Proposition}

\theoremstyle{definition}
\theoremstyle{definition}

\newcommand{\Z}{\mathbb{Z}}
\newcommand{\R}{\mathbb{R}}

\newcommand{\G}{\mathcal G}

\newcommand{\No}{\mathbf{N}}

\newcommand{\N}{\mathbb{N}}

\newcommand{\X}{X}

\newcommand{\eps}{\varepsilon}

\newcommand{\old}[1]{}
\renewcommand{\th}{\ensuremath{^{\text{th}}}\xspace}

\newcommand{\fm}{\phantom{-}}

\def\rcs $#1: #2 ${\expandafter\def\csname rcs#1\endcsname {#2}}
\rcs $Date: 2011/06/02 01:55:54 $

\begin{document}
\title{
Double-dimer pairings and skew Young diagrams
}
\author{\begin{tabular}{c}
\href{http://www.math.brown.edu/~rkenyon/}{Richard W. Kenyon}\\[-3pt]
\small Brown University\\[-4pt]
\small Providence, RI 02912, USA
\end{tabular} \and
\begin{tabular}{c}
\href{http://dbwilson.com}{David B. Wilson}\\[-3pt]
\small Microsoft Research\\[-4pt]
\small Redmond, WA 98052, USA
\end{tabular}
}
 \date{\dateline{July 12, 2010}{June 1, 2011}{June 14, 2011}}
\maketitle

\begin{abstract}
  We study the number of tilings of skew Young diagrams by ribbon
  tiles shaped like Dyck paths, in which the tiles are ``vertically
  decreasing''.  We use these quantities to compute pairing
  probabilities in the double-dimer model: Given a planar bipartite
  graph~$G$ with special vertices, called nodes, on the outer face,
  the double-dimer model is formed by the superposition of a uniformly
  random dimer configuration (perfect matching) of~$G$ together with a
  random dimer configuration of the graph formed from~$G$ by deleting
  the nodes.  The double-dimer configuration consists of loops,
  doubled edges, and chains that start and end at the boundary nodes.
  We are interested in how the chains connect the nodes.  An
  interesting special case is when the graph is $\varepsilon(\mathbb
  Z\times\mathbb N)$ and the nodes are at evenly spaced locations on
  the boundary~$\mathbb R$ as the grid spacing~$\varepsilon\to0$.
\end{abstract}

\section{Introduction}

\newcommand{\op}{\texttt{(}}
\newcommand{\cp}{\texttt{)}}

Among the combinatorial objects counted by Catalan numbers $C_n =
(2n)!/(n!(n+1)!)$ are balanced parentheses expressions (BPEs), Dyck
paths, and noncrossing pairings (see \cite[exercise~6.19]{MR1676282}).  A
word of length $2n$ in the symbols ``\op'' and ``\cp'' is said to be a
\textbf{balanced parentheses expression} if it can be reduced to the
empty word by successively removing subwords ``\op\cp''.
Equivalently, there are $n$ \op's and $n$ \cp's and the number of
\op's minus the number of \cp's in any initial segment is nonnegative
(see \cite[ex.~6.19(r)]{MR1676282}).  A \textbf{Dyck path} of length
$2n$ is a map $h:\{0,1,2,\dots,2n\}\to\{0,1,\dots\}$ with
$h(0)=h(2n)=0$ and $|h(i+1)-h(i)|=1$.  There is a bijection between
Dyck paths and BPEs defined by letting $h(i)$ be the number of \op's
minus the number of \cp's in the prefix of length $i$ (see
\cite[ex.~6.19(i)]{MR1676282}).  Dyck paths $h$ are also in bijective
correspondence with \textbf{noncrossing pairings} $\pi$, that is,
pairings of points $\{1,2,\dots,2n\}$ arranged counter-clockwise on a
circle, in which no two matched pairs cross.  The bijection is defined
as follows: the location of each up step is paired with the first
location at which the path returns to its previous height before the
up step (see \cite[ex.~6.19(n)]{MR1676282}).

These sets are also in bijection with ``confining'' subsets
$S\subseteq\{1,\dots,2n\}$: a subset $S$ is \textbf{confining} if it
has the same number of odd and even elements, and for any $i\in S^c =
\{1,\dots,2n\}\setminus S$, the set $S$ contains strictly more odds less than
$i$ than evens less than $i$.  (Confining sets are relevant to the
double-dimer model, as discussed in the next section, and there the
reason for this name will become clear.)  The bijection is as follows:
odd elements in $S$ and even elements in $S^c$ are replaced by \op;
even elements in $S$ and odd elements in $S^c$ are replaced by
\cp. The argument is left to the reader.

For any Dyck path $h$, let $S_h\subseteq\{1,\dots,2n\}$ be its
associated confining subset and let $\pi_h$ be its associated planar
pairing.

For example, when $n=3$, these bijections between the confining sets
$S$, BPEs, Dyck paths $h$, and planar pairings $\pi$ are summarized in
the following table.  The pairs in the pairing $\pi$ correspond to the
horizontal \textbf{chords} underneath the Dyck path which connect each
up step with its corresponding down step (shown in the diagrams).  The
table is arranged so that the Dyck paths $h$ are in lexicographic
order (equivalently the BPE's are in lexicographic order).
\begin{center}
\begin{tabular}{|c|c|c|c|c|}
\hline
confining set $S$ & BPE $P$ & Dyck path $h$ & pairing $\pi$ & diagram\\
\hline
$\{1,2,3,4,5,6\}$                     & $\op\cp\op\cp\op\cp$ & $0,1,0,1,0,1,0$ & $\,{}^1_2\!\mid\!{}^3_4\!\mid\!{}^5_6\,$ & \raisebox{-3pt}{\includegraphics{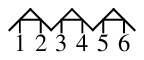}} \\
$\{1,2,3,\phantom{4,5,{}}6\}$           & $\op\cp\op\op\cp\cp$ & $0,1,0,1,2,1,0$ & $\,{}^1_2\!\mid\!{}^3_6\!\mid\!{}^5_4\,$ & \raisebox{-3pt}{\includegraphics{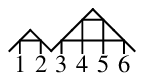}} \\
$\{1,\phantom{2,3,{}}4,5,6\}$           & $\op\op\cp\cp\op\cp$ & $0,1,2,1,0,1,0$ & $\,{}^1_4\!\mid\!{}^3_2\!\mid\!{}^5_6\,$ & \raisebox{-3pt}{\includegraphics{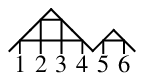}} \\
$\{1,\phantom{2,3,4,5,{}}6\}$           & $\op\op\cp\op\cp\cp$ & $0,1,2,1,2,1,0$ & $\,{}^1_6\!\mid\!{}^3_2\!\mid\!{}^5_4\,$ & \raisebox{-3pt}{\includegraphics{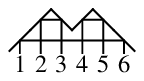}} \\
$\{1,\phantom{2,{}}3,4,\phantom{5,{}}6\}$ & $\op\op\op\cp\cp\cp$ & $0,1,2,3,2,1,0$ & $\,{}^1_6\!\mid\!{}^3_4\!\mid\!{}^5_2\,$ & \raisebox{-3pt}{\includegraphics{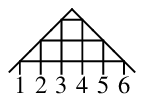}} \\
\hline
\end{tabular}
\end{center}

\newcommand{\rf}[1]{\begin{rotatebox}{90}{$#1$}\end{rotatebox}}
\newcommand{\n}{$-1\fm$}

\newcommand{\succp}{\stackrel{\op\cp}{\rightarrow}}
\newcommand{\precp}{\stackrel{\op\cp}{\leftarrow}}

We define a binary relation $\precp$ on BPEs (and therefore on Dyck
paths, etc.)  as follows.  We say $P_1\precp P_2$ if $P_1$ can be
obtained from $P_2$ by taking some of the matched pairs of parentheses
of $P_2$, and reversing each of them.  For example
$\op\cp\op\cp\precp\op\op\cp\cp$ by reversing the central pair of
matched parentheses. However
$\op\cp\op\op\cp\cp\,\not\!\precp\op\op\op\cp\cp\cp$.  The relation
$\precp$ is acyclic and reflexive, and its transitive closure is the well-known
partial order $\preceq$ on Dyck paths, where $h_1\preceq h_2$ if
$h_1(i)\leq h_2(i)$ for all $i$.  The lexicographic order used in
the above table is a linear extension of the partial order $\preceq$.

The following lemma helps motivate the binary relation $\precp$, and
we will use it in \sref{sec:ddimer} when we study the double-dimer
model.
\begin{lemma}\label{piS}
  Let $P_1,P_2$ be BPEs of equal length, $S_1$ the confining subset associated to
  $P_1$ and $\pi_2$ the pairing associated to $P_2$.  Then $P_1\precp
  P_2$ if and only if $\pi_2$ has no connection from $S_1$ to $S_1^c$.
\end{lemma}

\begin{proof}
  Suppose we reverse pairs $p_1,\dots,p_k$ of $P_2$ to make $P_1$.
  Each pair $p_i$ is a pair of $\pi_2$.  If $p_i=\{a,b\}$, then $a$
  and $b$ have opposite parity.  At one of these locations $h_2$ made
  an up step, and at the other $h_2$ made a down step, so either
  $a,b\in S_2$ or else $a,b\notin S_2$.  Reversing this pair, we have
  $a,b\notin S_1$ or $a,b\in S_1$, so the pairs $p_1,\dots,p_k$ of
  $\pi_2$ do not match $S_1$ to $S_1^c$.  For any pair $\{a,b\}$ of
  $\pi_2$ that is not reversed, either $a,b\in S_2$ or else $a,b\notin
  S_2$, in which case $a,b\in S_1$ or else $a,b\notin S_1$, so these
  pairs do not match $S_1$ to $S_1^c$ either.

  Conversely, suppose $\pi_2$ has no connections from $S_1$ to
  $S_1^c$.  Consider a pair $\{a,b\}$ of $\pi_2$ where $a<b$ (so $a$
  and $b$ have opposite parity and $a$ is an up step of $h_2$ and $b$
  is a down step of $h_2$).  If $a$ is an up step of $h_1$, then since
  $\{a,b\}$ does not connect $S_1$ to $S_1^c$, it must be that $b$ is
  a down step of $h_1$.  Likewise, if $a$ is a down step of $h_1$,
  then $b$ is an up step of $h_1$.  Thus if we reverse the parentheses
  in BPE $P_2$ for all such pairs $\{a,b\}$ for which $a$ is a down
  step in $h_1$, then the result is $P_1$.
\end{proof}

\subsection{The incidence matrix \texorpdfstring{$M$}{M} and its inverse\texorpdfstring{ $M^{-1}$}{}}

Associated to the binary relation $\precp$ is its ``incidence''
matrix~$M$, with $M_{P_1,P_2}=\delta_{\{P_1\precp P_2\}}$.  It is
convenient to order the rows and columns according the lexicographic
order on BPEs.  Since this order is a linear extension of $\preceq$,
which in turn is the transitive closure of $\precp$, the matrix $M$
will be upper triangular when written in this way.  For example, when
$n=3$, $M$ is
$$
%\begin{gathered}[b]
\begin{array}{rrrrrrr}
               &                                           &
\rf{\op\cp\op\cp\op\cp \ \,{}^1_2\!\mid\!{}^3_4\!\mid\!{}^5_6\,} &
\rf{\op\cp\op\op\cp\cp \ \,{}^1_2\!\mid\!{}^3_6\!\mid\!{}^5_4\,} &
\rf{\op\op\cp\cp\op\cp \ \,{}^1_4\!\mid\!{}^3_2\!\mid\!{}^5_6\,} &
\rf{\op\op\cp\op\cp\cp \ \,{}^1_6\!\mid\!{}^3_2\!\mid\!{}^5_4\,} &
\rf{\op\op\op\cp\cp\cp \ \,{}^1_6\!\mid\!{}^3_4\!\mid\!{}^5_2\,}
\\
\{1,2,3,4,5,6\}                       & \op\cp\op\cp\op\cp & 1& 1& 1& 1& 1 \\[2pt]
\{1,2,3,\phantom{4,5,{}}6\}           & \op\cp\op\op\cp\cp & 0& 1& 0& 1& 0 \\[2pt]
\{1,\phantom{2,3,{}}4,5,6\}           & \op\op\cp\cp\op\cp & 0& 0& 1& 1& 0 \\[2pt]
\{1,\phantom{2,3,4,5,{}}6\}           & \op\op\cp\op\cp\cp & 0& 0& 0& 1& 1 \\[2pt]
\{1,\phantom{2,{}}3,4,\phantom{5,{}}6\}&\op\op\op\cp\cp\cp & 0& 0& 0& 0& 1 \\[2pt]
%\end{matrix}
\end{array}
%\end{gathered}
$$
In addition to labeling the rows and columns of $M$ by the BPEs, we have also
labeled the rows of $M$ according the corresponding confining subsets, and the
columns by the corresponding noncrossing pairings, since this is how we shall use
the matrices in \sref{sec:ddimer}.

Since the matrix $M$ is upper triangular with $1$'s on the diagonal,
it is invertible and $M^{-1}$ is upper triangular with integer
entries.  The matrix $M^{-1}$ is analogous to the M\"obius function of
a partial order, except that the binary relation associated with $M$
is not transitive.  For $n=3$, $M^{-1}$ is
$$
\begin{array}{rrrrrrr}
%\begin{matrix}
 &&
\rf{\op\cp\op\cp\op\cp \ \{1,2,3,4,5,6\}} &
\rf{\op\cp\op\op\cp\cp \ \{1,2,3,\phantom{4,5,{}}6\}} &
\rf{\op\op\cp\cp\op\cp \ \{1,\phantom{2,3,{}}4,5,6\}} &
\rf{\op\op\cp\op\cp\cp \ \{1,\phantom{2,3,4,5,{}}6\}} &
\rf{\op\op\op\cp\cp\cp \ \{1,\phantom{2,{}}3,4,\phantom{5,{}}6\}}
\\
\,{}^1_2\!\mid\!{}^3_4\!\mid\!{}^5_6\, & \op\cp\op\cp\op\cp & 1& -1& -1& 1& -2 \\[2pt]
\,{}^1_2\!\mid\!{}^3_6\!\mid\!{}^5_4\, & \op\cp\op\op\cp\cp & 0& 1& 0& -1& 1   \\[2pt]
\,{}^1_4\!\mid\!{}^3_2\!\mid\!{}^5_6\, & \op\op\cp\cp\op\cp & 0& 0& 1& -1& 1   \\[2pt]
\,{}^1_6\!\mid\!{}^3_2\!\mid\!{}^5_4\, & \op\op\cp\op\cp\cp & 0& 0& 0& 1& -1   \\[2pt]
\,{}^1_6\!\mid\!{}^3_4\!\mid\!{}^5_2\, & \op\op\op\cp\cp\cp & 0& 0& 0& 0& 1    \\[2pt]
%\end{matrix}
\end{array}
$$

\subsection{Skew Young diagrams}

We may associate with each Dyck path an integer partition, or
equivalently, a Young diagram, which is given by the set of boxes
which may be placed above the Dyck path.  Dyck paths of different
lengths may be associated to the same Young diagram, but distinct Dyck
paths of the same length will be associated to distinct Young diagrams.

The matrices $M$ and $M^{-1}$ can be interpreted as being indexed by
pairs of Dyck paths, which in turn correspond to pairs of Young
diagrams.  If $M_{\lambda,\mu}$ is
non-zero, then $\mu$ is larger than $\lambda$ as a path, or
equivalently $\mu$ is smaller than $\lambda$ as a partition
($\mu\subseteq\lambda$).  Likewise, since $M^{-1}$ is also upper triangular,
if $M^{-1}_{\lambda,\mu}$ is nonzero then $\mu\subseteq\lambda$.
Each matrix entry can be associated with the
skew Young diagram $\lambda/\mu$.  Different matrix entries can
correspond to the same skew Young diagram, but as the next two (easy)
lemmas show, when this happens, the matrix entries are equal.

\begin{lemma}
  Suppose that $\mu_1\subseteq\lambda_1$ and
  $\mu_2\subseteq\lambda_2$, and the skew shapes $\lambda_1/\mu_1$ and
  $\lambda_2/\mu_2$ are equivalent in the sense that $\lambda_1/\mu_1$
  may be translated to obtain $\lambda_2/\mu_2$.  Then
  $M_{\lambda_1,\mu_1} = M_{\lambda_2,\mu_2}$.
\end{lemma}
\begin{proof}
  Suppose $M_{\lambda_1,\mu_1}=1$.  Then $\lambda_1$ may be obtained
  from $\mu_1$ by ``pushing down'' on some of $\mu_1$'s chords.  Each
  such chord of $\mu_1$ must lie within a connected component of the skew shape
  $\lambda_1/\mu_1$.  The result is now easy.
\end{proof}

Therefore we may define $M_{\lambda/\mu}$ to be $M_{\lambda,\mu}$.

\begin{lemma}
  $M^{-1}_{\lambda,\mu}$ is determined by the translation equivalence
  class of the skew shape $\lambda/\mu$.
\end{lemma}
\begin{proof}
  We prove this by induction on $|\lambda/\mu|$.  We have
\begin{align*}
   \sum_{\rho:\mu\subseteq\rho\subseteq\lambda} M_{\lambda/\rho} M^{-1}_{\rho,\mu} &= \delta_{\mu,\lambda}, \\
\intertext{and isolating the $\rho=\lambda$ term,}
   M^{-1}_{\lambda,\mu} &= \delta_{\mu,\lambda} - \sum_{\rho:\mu\subseteq\rho\subsetneq\lambda} M_{\lambda/\rho} M^{-1}_{\rho/\mu},
\end{align*}
which only depends on $\lambda/\mu$.
\end{proof}

\begin{figure}[htbp]
\def\sz{0.9}
\psfrag{q1}[cc][cc][\sz][0]{$-1$}
\psfrag{q1^2}[cc][cc][\sz][0]{$1$}
\psfrag{q1^4}[cc][cc][\sz][0]{$1$}
\psfrag{q1^3}[cc][cc][\sz][0]{$-1$}
\psfrag{q1^5}[cc][cc][\sz][0]{$-1$}
\psfrag{q3q1+q1^4}[cc][cc][\sz][0]{$2$}
\psfrag{q3+q1^3}[cc][cc][\sz][0]{$-2$}
\psfrag{q3q1^4+q1^7}[cc][cc][\sz][0]{$-2$}
\psfrag{q3q1^3+q1^6}[cc][cc][\sz][0]{$2$}
\psfrag{q3q1^2+q1^5}[cc][cc][\sz][0]{$-2$}
\psfrag{q3q1+q1^4}[cc][cc][\sz][0]{$2$}
\psfrag{q5+q3q1^2+q1^5}[cc][cc][\sz][0]{$-3$}
\psfrag{q5+2q3q1^2+q1^5}[cc][cc][\sz][0]{$-4$}
\psfrag{q3^2q1^2+q3q1^5+q1^8}[cc][cc][\sz][0]{$3$}
\includegraphics[width=\textwidth]{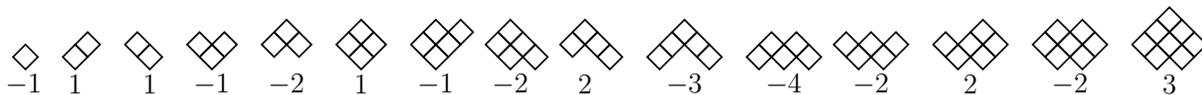}
\caption{The values of $M^{-1}_{\lambda/\mu}$ for the first few
  skew shapes $\lambda/\mu$.
  \label{skews}}
\end{figure}

Thus the expression $M^{-1}_{\lambda/\mu}$ is well-defined.  We give
in \fref{skews} $M^{-1}_{\lambda/\mu}$ for the first few
connected skew shapes.  The next lemma shows that for
disconnected $\lambda/\mu$, $M^{-1}_{\lambda/\mu}$ is a product
of the $M^{-1}$'s for the connected components of $\lambda/\mu$.

\begin{lemma} \label{m-1-factor}
  Suppose $\mu\subseteq\lambda$ and $\lambda/\mu$ has $k$ connected
  components $\lambda_1/\mu_1,\dots,\lambda_k/\mu_k$.  Then
  $M_{\lambda/\mu} = M_{\lambda_1/\mu_1} \times\dots\times
  M_{\lambda_k/\mu_k}$ and
  $M^{-1}_{\lambda/\mu} = M^{-1}_{\lambda_1/\mu_1} \times\dots\times
  M^{-1}_{\lambda_k/\mu_k}$.
\end{lemma}
\begin{proof}
  If we can reverse parentheses in $\mu$'s BPE to obtain $\lambda$'s
  BPE, then in the Dyck path representation, the chords connecting the
  Dyck path steps corresponding to the parentheses to be reversed will
  lie within the region $\lambda/\mu$.  The multiplicative property
  for $M_{\lambda/\mu}$ follows.  For $M^{-1}_{\lambda/\mu}$, the multiplicative
  property follows from induction and the following equation
\begin{align*}
   \sum_{\substack{\rho_1:\mu_1\subseteq\rho_1\subseteq\lambda_1 \\ \vdots \\ \rho_k:\mu_k\subseteq\rho_k\subseteq\lambda_k}} \underbrace{M_{\lambda_1/\rho_1} \cdots M_{\lambda_k/\rho_k}}_{M_{\lambda/\rho}} M^{-1}_{\rho_1/\mu_1} \cdots M^{-1}_{\rho_k/\mu_k} &= \delta_{\mu_1,\lambda_1} \cdots \delta_{\mu_k,\lambda_k} = \delta_{\mu,\lambda}. \qedhere
\end{align*}
\end{proof}

\begin{figure}[htbp]
\centerline{\includegraphics{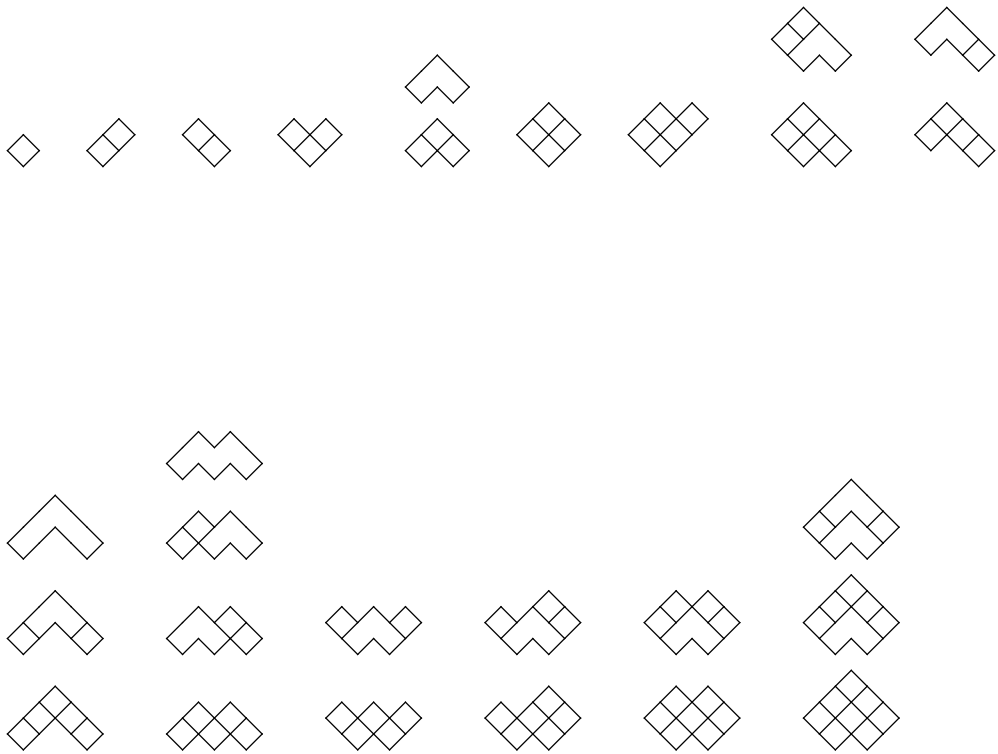}}
\caption{Cover-inclusive Dyck tilings of the first few skew shapes.}
\label{skew-small}
\end{figure}

\begin{figure}
\centerline{\includegraphics{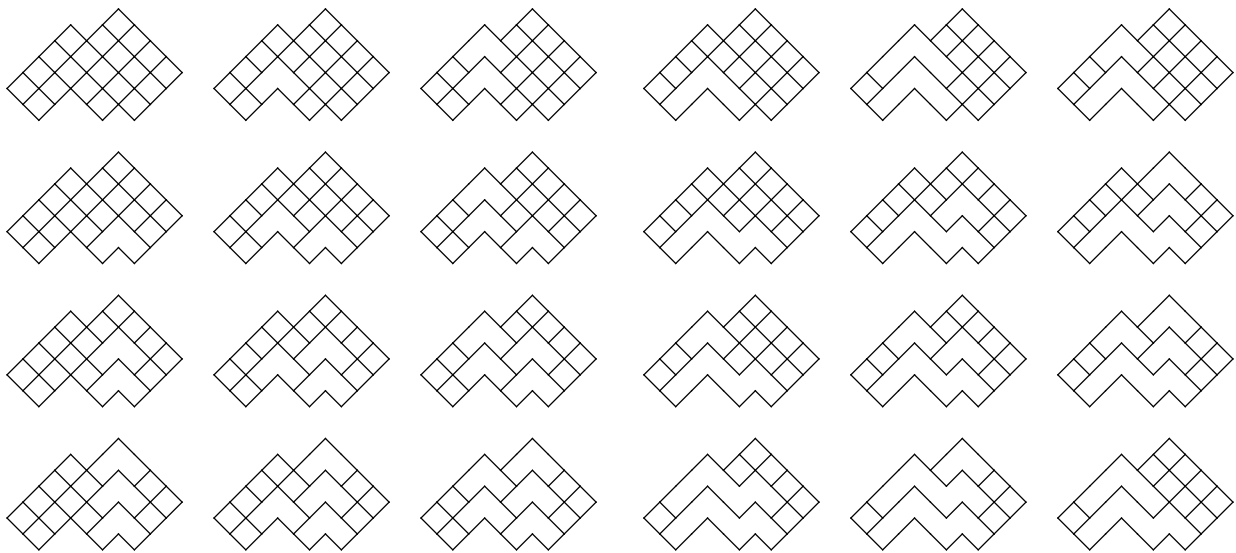}}
\caption{Cover-inclusive Dyck tilings of a larger skew shape.}
\label{skew-bigger}
\end{figure}

The following theorem gives a formula for $M^{-1}_{\lambda/\mu}$.
The sign is given by the parity of the area of the skew shape.  The
absolute value is the number of certain tilings, as indicated in
Figures~\ref{skew-small} and~\ref{skew-bigger}:
Each tile is essentially an expanded version of a Dyck path,
where each point in a Dyck path is replaced with a box, so we call it
a \textbf{Dyck tile}.  Dyck tiles are ribbon tiles in which the start
box and end box are at the same height, and no box within the tile is
below them.  A tiling of the skew Young diagram by Dyck tiles is a
\textbf{Dyck tiling}.  We say that one Dyck tile covers another Dyck
tile if the first tile has at least one box whose center lies straight above
the center of a box in the second tile.  We say that a Dyck tiling is
\textbf{cover-inclusive} if for each pair of its tiles, if the first tile
covers the second tile, then the horizontal extent of the first tile
is included as a subset within the horizontal extent of the second tile.
We shall prove the following theorem:
\begin{theorem} \label{thm:tilings}
$ M^{-1}_{\lambda/\mu} = (-1)^{|\lambda/\mu|} \times |\{\text{cover-inclusive Dyck tilings of $\lambda/\mu$}\}|$.
\end{theorem}

To prove this, we start with a recursive formula for computing
$M^{-1}_{\lambda/\mu}$.
Using the fact $M^{-1} M = I$, we get
\begin{align*}
   \sum_{\rho} M^{-1}_{\lambda,\rho} M_{\rho,\mu} &= \delta_{\mu,\lambda} \\
   \sum_{\rho:\mu\subseteq\rho\subseteq\lambda} M^{-1}_{\lambda/\rho} M_{\rho/\mu} &= \delta_{\mu,\lambda}.
\end{align*}
Let us consider the chords of $\mu$ that start and end
within the region of the skew shape $\lambda/\mu$.  Rather than
summing over all $\rho$'s for which
$\mu\subseteq\rho\subseteq\lambda$, because of the $M_{\rho/\mu}$ factor,
we can instead sum over all
$\rho$'s which may be obtained from $\mu$ by pushing down on some of
the chords which start and end within the region of $\lambda/\mu$.
\begin{align}
   \sum_{\rho\subseteq\lambda:\text{$\rho$ obtained by pushing down chords of $\mu$}} M^{-1}_{\lambda/\rho} &= \delta_{\mu,\lambda}.
\end{align}
This is the ``downward recurrence''.
\begin{figure}[b]
\psfrag{(a)}{(a)}
\psfrag{(b)}{(b)}
\psfrag{(c)}{(c)}
\psfrag{(d)}{(d)}
\psfrag{(e)}{(e)}
\includegraphics{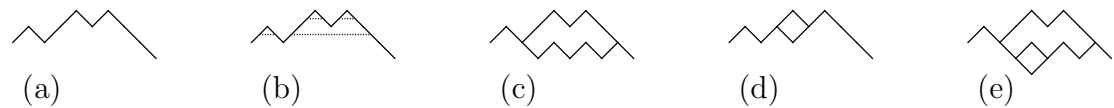}
\caption{Illustration of pushing chords down and the resulting Dyck tiles: (a) an example upper boundary of a skew shape $\lambda/\mu$ (the upper boundary need not be a Dyck path, but can be extended to a Dyck path), (b) the chords of this upper boundary, corresponding to pairs of parentheses that can be reversed, (c) the boundary together with the result of pushing down the long chord, (d) the upper boundary together with the result of pushing down the middle short chord, (e) if multiple chords are pushed down, the outermost chords are pushed down first.  Pushing down chords is equivalent to laying down Dyck tiles along the upper bondary.
  \label{push}}
\end{figure}
We get a different recurrence, the ``upward recurrence'' from $MM^{-1}  = I$:
\begin{align}
   \sum_{\rho} M_{\lambda,\rho} M^{-1}_{\rho,\mu} &= \delta_{\mu,\lambda} \notag\\
   \sum_{\rho\supseteq\mu:\text{$\lambda$ obtained by pushing down chords of $\rho$}} M^{-1}_{\rho/\mu} &= \delta_{\mu,\lambda}.
\end{align}
Let us rewrite the downward recurrence for $\lambda/\mu\neq\varnothing$:
\begin{align}
 M^{-1}_{\lambda/\mu} =  \sum_{\text{nonempty sets $S$ of chords of $\mu$}} -M^{-1}_{\lambda/(\text{$\mu$ with $S$ pushed down})}.
\end{align}
Rather than restrict to sets of chords of~$\mu$ that can be pushed
down to obtain a path~$\rho$ above~$\lambda$, it is convenient to sum
over all possible (nonempty) sets of chords, with the understanding
that $M^{-1}_{\lambda/\rho}=0$ if $\rho\not\subseteq\lambda$.

Every time we push a chord of $\mu$, we can interpret that as laying
down a Dyck tile along the upper boundary (adjacent to $\mu$) of the skew shape
$\lambda/\mu$ (see \fref{push}).  If multiple chords are pushed, then multiple Dyck
tiles are laid down, where we follow the convention that longer
chords are laid down first (so the shorter ones are below the longer ones).
If we expand the recursive formula for
$M^{-1}_{\lambda/\mu}$, then each term corresponds to a Dyck tiling
of $\lambda/\mu$.  It is convenient to let each Dyck tile have
weight~$-1$.  Since each tile contains an odd number of boxes, the
parity of the $-1$ factors in a tiling is just the parity of
$|\lambda/\mu|$.  We further rewrite the downward recurrence as
\begin{align}
(-1)^{|\lambda/\mu|} M^{-1}_{\lambda/\mu} =  \sum_{\substack{\text{nonempty sets $S$ of Dyck tiles}\\ \text{that can be placed along upper edge of $\mu$}}} (-1)^{1+|S|} \big[(-1)^{|\lambda/(\mu\downarrow S)|}M^{-1}_{\lambda/(\mu\downarrow S)}\big] \label{recurrenceMinv},
\end{align}
where $(\mu\downarrow S)$ denotes $\mu$ with $S$ pushed down.

\begin{figure}[htbp]
\def\sz{1.0}
\psfrag{+}[cc][cc][\sz][0]{$+$}
\psfrag{-}[cc][cc][\sz][0]{$-$}
\psfrag{=}[cc][cc][\sz][0]{$=$}
\psfrag{3}[cc][cc][\sz][0]{\vspace*{3pt}$\overbrace{\hspace{1.15in}}$}
\psfrag{4}[cc][cc][\sz][0]{\vspace*{3pt}$\overbrace{\hspace{1.53in}}$}
\psfrag{9}[cc][cc][\sz][0]{\vspace*{3pt}$\overbrace{\hspace{3.6in}}$}
\psfrag{A}[cc][cc][\sz][0]{\vspace*{3pt}$\overbrace{\hspace{4.0in}}$}
\includegraphics[width=\textwidth]{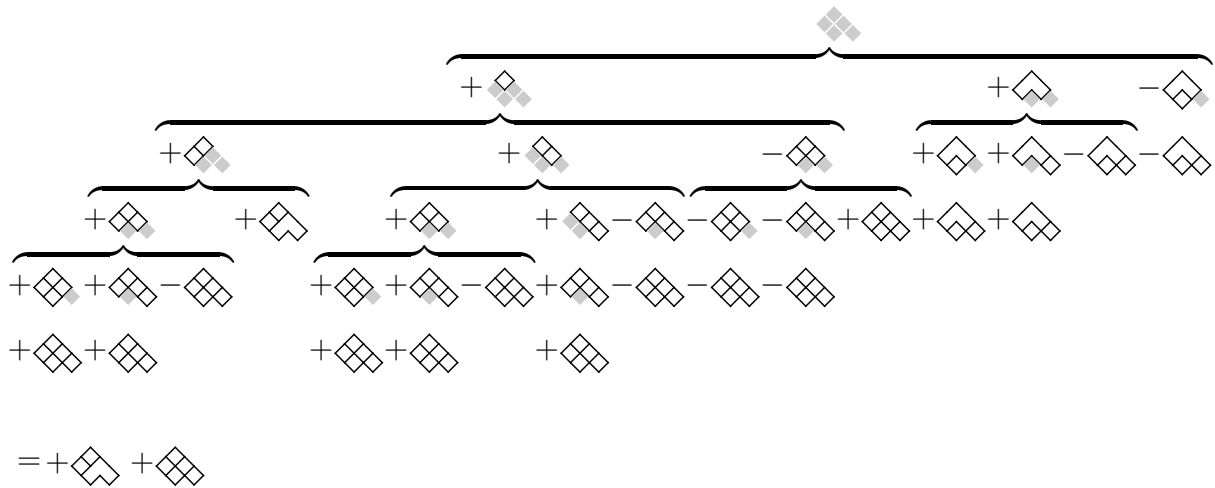}
\caption{Example of using the recursive definition of $f_{\lambda/\mu}$.  There is a lot of cancelation, and the result is the sum of cover-inclusive Dyck tilings.
  \label{skew-recursive}}
\end{figure}

Let us define $f_{\lambda/\mu}$ to be the formal linear combination of
Dyck tilings of the skew shape $\lambda/\mu$ defined recursively by
\begin{align}\label{recurrencef}
f_{\lambda/\mu} =  \sum_{\substack{\text{nonempty sets $S$ of Dyck tiles}\\ \text{that can be placed along upper edge of $\mu$}}} (-1)^{1+|S|} \left(\parbox{2.4in}{tiles of $S$ placed on top of $f_{\lambda/(\mu\downarrow S)}$,\\ with longer tiles higher up}\right),
\end{align}
with the base cases $f_{\lambda/\mu}=1$ if $\lambda/\mu=\varnothing$
and $f_{\lambda/\mu}=0$ if $\mu\not\subseteq\lambda$.  Comparing the
recurrences \eqref{recurrenceMinv} and \eqref{recurrencef}, we see
that if each Dyck tiling is replaced with~$1$, then $f_{\lambda/\mu}$
simplifies to $(-1)^{|\lambda/\mu|} M^{-1}_{\lambda/\mu}$.  In view of
this, we see that \tref{thm:tilings} is a corollary of
\tref{thm:tilings-formal}:

\begin{theorem} \label{thm:tilings-formal}
  With $f_{\lambda/\mu}$ defined in \eqref{recurrencef}, $f_{\lambda/\mu}$ is in fact a
  linear combination of just the cover-inclusive Dyck tilings, with a
  coefficient of~$1$ for each such tiling:
 $$ f_{\lambda/\mu} = \sum_{\text{cover-inclusive Dyck tilings $\mathcal T$ of $\lambda/\mu$}} \mathcal T.$$
(See \fref{skew-recursive} for an example.)
\end{theorem}

To prove \tref{thm:tilings-formal}, we use the following two lemmas.

\begin{lemma} \label{lem:top-tile}
  Any Dyck tiling of a skew Young diagram $\lambda/\mu
  \neq\varnothing$ contains a tile $T$ along the upper boundary of
  $\lambda/\mu$ such that $(\lambda/\mu)\setminus T$ is a skew Young
  diagram.
\end{lemma}
\begin{proof}
  Let $T_1$ be the tile of the Dyck tiling which contains the
  left-most square along the upper boundary of $\lambda/\mu$.  Suppose
  tile $T_n$ borders the upper boundary of $\lambda/\mu$ and has no
  tile above the upper-left edge of its leftmost square (e.g., $T_1$).
  Either $(\lambda/\mu)\setminus T_n$ is a skew Young diagram, or else
  we may define $T_{n+1}$ to be the tile containing the leftmost
  square square that borders the upper boundary of $T_n$.  Tile
  $T_{n+1}$ borders the upper boundary of $\lambda/\mu$, has no tile
  above the upper-left edge of its leftmost square, and its leftmost
  square is to the right of the leftmost square of $T_n$.  Since there
  are finitely many tiles in the Dyck tiling, the lemma follows.
\end{proof}

Given a Dyck tiling of a skew Young diagram $\lambda/\mu$ with $k$
tiles, define a valid labeling to be an assignment of the numbers
$1,\dots,k$ to the tiles such that for each $j\leq k$, tiles
$1,\dots,j$ form a skew Young diagram with tile $j$ along its upper
boundary.  (By \lref{lem:top-tile} such labelings exist.)  Given
two tiles $T_1$ and $T_2$ in the Dyck tiling, say that $T_1\prec T_2$
if $T_1$'s label is smaller than $T_2$'s label in all such labelings.
Then $\prec$ is a partial order on the tiles of the Dyck tiling.

\begin{lemma} \label{lem:cover}
  Suppose a Dyck tiling of a skew Young diagram $\lambda/\mu$ contains
  no pair of tiles $T_1$ and $T_2$ for which (1) $T_2$ is above $T_1$
  (2) $T_1$ and $T_2$ have no tiles between them, and (3) the
  horizontal extent of $T_1$ is a proper subset of the extent of
  $T_2$.  Then the tiling is cover-inclusive.
\end{lemma}
\begin{window}[2,r,\includegraphics{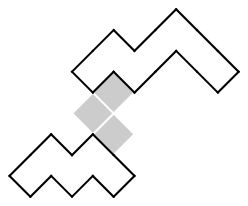},{}]
  \noindent\textit{Proof.}
  We prove the contrapositive.  If a Dyck tiling is not cover
  inclusive, then there is a pair of tiles $A$ and $B$, with $B$ above
  $A$, for which the horizontal extent of $B$ is not a subset of the
  extent of $A$.  If the set of squares above $A$ and below $B$ (shown
  in gray in the figure) is
  nonempty, then let $C$ be any tile containing a
  square between $A$ and $B$.  If the horizontal extent of $B$ is not
  a subset of the extent of $C$, then we may consider instead the pair
  of tiles $C$ and $B$, and if the extent of $B$ is a subset of the
  extent of $C$, then we may consider instead the pair of tiles $A$
  and $C$.  For this new pair of tiles, the extent of the upper tile
  is not a subset of the extent of the lower tile, and the interval with
  respect to the partial order $\prec$ between the new pair of tiles
  is strictly smaller than the interval w.r.t.\ $\prec$ between $A$
  and $B$.  Thus by induction we may assume that tiles $A$ and $B$
  have no squares between them.
\end{window}

\begin{window}[0,r,\includegraphics{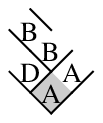},{}]
  If the extent of $A$ is a (proper) subset of the extent of $B$,
  then we may take $T_1=A$ and $T_2=B$.
  Otherwise, exactly one of the endpoints of $A$ lies within the extent of $B$;
  suppose without loss of generality it is the left, and let $L$
  denote the leftmost box of $A$ (shown in gray in the figure).
  Since $B$ does not cover $A$, tile
  $A$ contains the box immediately up and right of $L$.  Since $A$ and
  $B$ have no tiles between them, the box above $L$ must be part of
  tile $B$.  The right endpoint of tile $B$ lies above $A$, and since
  $B$ is a Dyck tile, the box to the immediate upper left of $L$
  cannot be part of tile $B$, say that it is part of tile $D$. This
  box is the right endpoint of tile $D$, and since $D$ is a Dyck tile
  and $B$ is a Dyck tile, the left endpoint of $B$ is to the left of
  $D$'s left endpoint (otherwise tile $D$ would have a square lower than
  its endpoints).  Thus, the horizontal extent of $D$ is a proper
  subset of the horizontal extent of $B$, so we take $T_1=D$ and $T_2=B$.
\hfil\qed
\end{window}
\vspace{2pt}

\begin{proof}[Proof of \tref{thm:tilings-formal}]
We prove the theorem by induction on $|\lambda/\mu|$.  We have equality when
$\lambda/\mu=\varnothing$.  Otherwise, there is some set $S$ of
possible Dyck tiles that may be placed along the upper boundary of
$\lambda/\mu$, which correspond to pushing a single chord of $\mu$ down.  If distinct
tiles $T_1,T_2\in S$ overlap, then one tile is a subset of the other,
say $T_1\subseteq T_2$ (this is because two chords of $\mu$ cannot have 
interleaved endpoints).  

Let us evaluate $f_{\lambda/\mu}$ restricted to tilings with $T_2$ at
the top edge and $T_1$ directly under it.  Either $T_1$ and $T_2$ are
pushed down in the same step of the recurrence (\ref{recurrencef}) or
$T_1$ is pushed down after $T_2$.

The subsets $S$ for the first step of the recurrence in which $T_2$ is
present at the top may be paired off with one another so that the
symmetric difference of each pair is $\{T_1\}$.  (Here we are not
assuming that $\mu\downarrow S$ lies above $\lambda$; when it does not
lie above $\lambda$ we have $f_{\lambda/(\mu\downarrow S)} = 0$.)
Let $A,A\cup\{T_1\}$ be such a pair.  Comparing
$f_{\lambda/(\mu\downarrow (A\cup\{T_1\}))}$ with
$f_{\lambda/(\mu\downarrow A)}$ restricted to tilings with $T_2$ at
the top edge and $T_1$ directly under it, by the induction hypothesis
they are exactly the same, but they have opposite signs in the
recursive formula for $f_{\lambda/\mu}$.  Thus $f_{\lambda/\mu}$ has
no tilings in which a tile $T_1$ is directly covered by a strictly
longer tile $T_2$.
Now using \lref{lem:cover},
we conclude that $f_{\lambda/\mu}$ contains only cover-inclusive Dyck tilings.

For a given cover-inclusive Dyck tiling $\mathcal T$ of skew shape
$\mu/\lambda$, let us find its coefficient in $f_{\lambda/\mu}$, which
we denote $[\mathcal T]f_{\lambda/\mu}$.  Let $S$ be the set of tiles along the
top edge of $\mathcal T$ that can be pushed down in the first step
of the recurrence.  (By \lref{lem:top-tile}, $S\neq\varnothing$.)
If any tile not in $S$ is pushed down, then the result cannot be extended to $\mathcal T$, so we have
  $$[\mathcal T]f_{\lambda/\mu} = \sum_{\substack{S'\subseteq S\\S'\neq\varnothing}} (-1)^{1+|S'|} [\mathcal T] f_{\lambda/(\mu\downarrow S')} = \sum_{\substack{S'\subseteq S\\S'\neq\varnothing}} (-1)^{1+|S'|} = 1,$$
which completes the induction in the proof of \tref{thm:tilings-formal}.
\end{proof}

Let us define $f_{\lambda/\mu}(q)$ to be the polynomial obtained by giving each tile weight~$q$.
Then 
\begin{equation}
  M^{-1}_{\lambda/\mu} = f_{\lambda/\mu}(-1).
\end{equation}

Given the cover-inclusive Dyck tiling characterization, the next several propositions are straightforward to verify. Recall the $q$-analogue notation
\begin{eqnarray*}
n_q&=&1+q+\dots+q^{n-1}\\
n!_q&=&n_q(n-1)_q\dots 1_q\\
\binom{a}{b}_q &=& \frac{a!_q}{b!_q(a-b)!_q}.
\end{eqnarray*}

\begin{proposition}
  If the lower boundary of the skew shape $\lambda/\mu$ is minimal ($V$-shaped) then $f_{\lambda/\mu}(q) = q^{|\lambda/\mu|}$.
$$
\raisebox{-24pt}{
\includegraphics[scale=0.7]{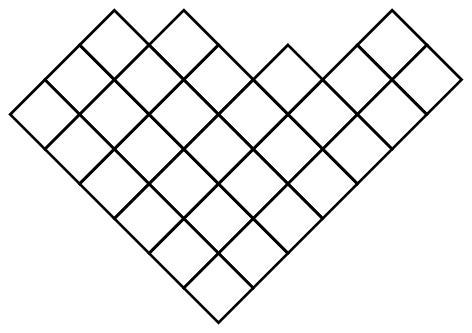}}
\Rightarrow
q^{|\lambda/\mu|}
$$
\end{proposition}
\begin{proof}[Proof sketch]
  There is only one cover-inclusive Dyck tiling; in it each tile has size~1.
\end{proof}

\begin{proposition}
  If $\lambda/\mu$ is $\Lambda$-shaped, then $f_{\lambda/\mu}(q)$
  is $q^{|\lambda/\mu|}$ times a $q^{-2}$-analogue of a
  binomial coefficient as illustrated in the
  following example:
$$
\raisebox{-24pt}{
\def\sz{0.6}
\psfrag{a=4}[cc][cc][\sz][0]{$a=4$}
\psfrag{b=2}[cc][cc][\sz][0]{$b=2$\ }
\psfrag{c=3}[cc][cc][\sz][0]{\ $c=3$}
\psfrag{d=3}[cc][cc][\sz][0]{$d=3$}
\includegraphics[scale=0.7]{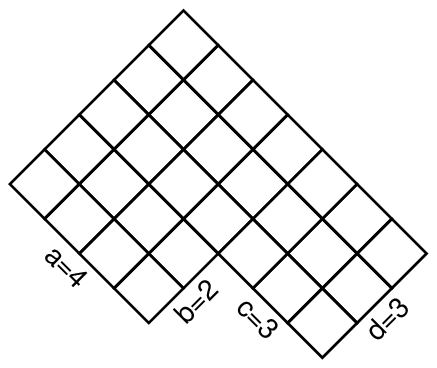}}
\Rightarrow
q^{|\lambda/\mu|}
\binom{\min(a,d)+\min(b,c)}{\min(a,d)}_{q^{-2}}
$$
\end{proposition}
\begin{proof}[Proof sketch]
  In any cover-inclusive Dyck tiling, the tiles with size larger
  than~1 are all $\Lambda$-shaped and centered at the peak on the
  lower boundary, and their sizes decrease (by even numbers) when going up.
\end{proof}

\begin{proposition}[First row of $M^{-1}$] \label{row1} If $\lambda$
  is the zigzag path, then $f_{\lambda/\mu}(q)$ is
  $q^{|\lambda/\mu|}$ times a product of $q^{-2}$-analogues of
  heights, as illustrated in the following example:
\begin{gather*}
\includegraphics[scale=0.7]{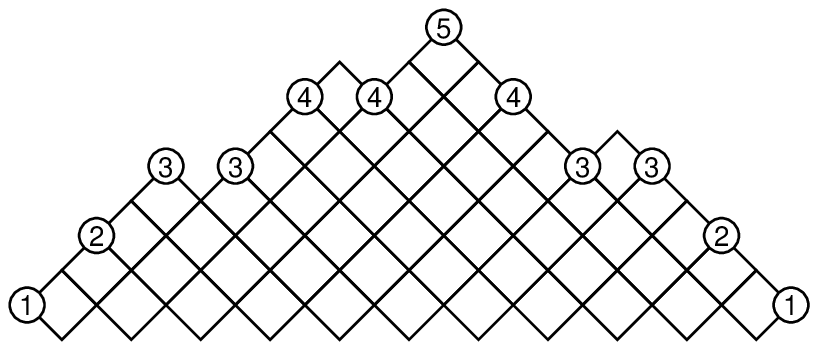}
\\\Rightarrow q^{52} \times
1_{q^{-2}}\times 2_{q^{-2}}\times 3_{q^{-2}}\times 3_{q^{-2}} \times 4_{q^{-2}}\times 4_{q^{-2}}\times 5_{q^{-2}}\times 4_{q^{-2}} \times 3_{q^{-2}}\times 3_{q^{-2}}\times 2_{q^{-2}}\times 1_{q^{-2}}
\end{gather*}
\end{proposition}
%\begin{proof}[Proof sketch]
\begin{window}[0,r,\includegraphics[scale=0.7]{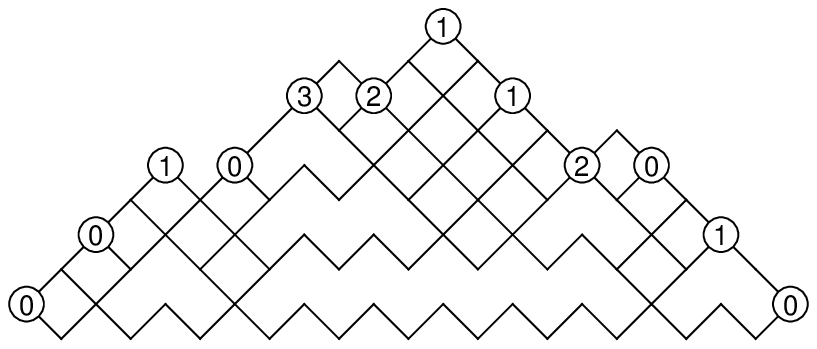}\\,{}]
  \noindent\textit{Proof sketch.}
  In any cover-inclusive Dyck tiling of such regions, each tile is a
  zigzag shape which starts at a square with even parity (when the
  lower-left-most square has even parity).  Each location with a
  circled number $h$ specifies a number between $0$ and $h-1$,
  which is the number of squares directly below it that get glued
  to their lower neighbors to be part of the same Dyck tile,
  as shown in the figure.
  These tower heights are independent of one
  another and determine the tiling. \hfill $\square$
\end{window}
%\end{proof}
\newpage

\begin{proposition}\label{snake}
  If $\lambda/\mu$ is a width-one strip, then $M^{-1}_{\lambda/\mu}$
  is up to sign given by a nested sequence of products and $+1$'s,
  with a $+1$ for each chord and a times for each minimum, as
  indicated in the figure.
\begin{gather*}
\psfrag{*}[cc][cc][0.5][0]{$)\times($}
\psfrag{+1}[cc][cc][0.5][0]{$+1$}
\psfrag{1}[cc][cc][0.5][0]{$1$}
\psfrag{op}[cc][cc][0.5][0]{$($}
\psfrag{cp}[cc][cc][0.5][0]{$)$}
\psfrag{cp*op}[cc][cc][0.5][0]{$)\times($}
\includegraphics[scale=0.7]{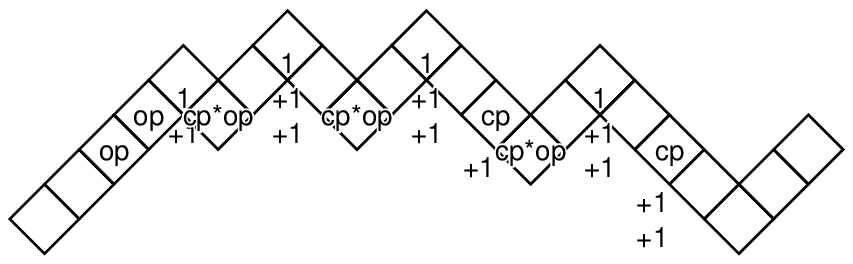}
\\
\Rightarrow
\pm[((1+1)\times(1+1+1)\times(1+1+1)+1)\times(1+1+1)+1+1]
\end{gather*}
\end{proposition}
\begin{proof}[Proof sketch]
  Each $+1$ term at a chord (or local maximum) corresponds to the
  endpoints of the chord belonging to the same Dyck tile (which
  extends no lower than the chord).  Each $\times$ at a local minimum
  corresponds to making independent choices to the left and right of
  the local minimum.
\end{proof}

We also observed some formulas (Conjectures~\ref{rowsum}
and~\ref{columnsum}) for which we do not have proofs.

\begin{conjecture} \label{rowsum}
  For any row~$\lambda$ of $M^{-1}$ (of order $n$), the absolute values of the
  entries add up to a divisor of $n!$. If the lengths of the
  chords are measured as $1,2,3,\dots$, then the $\lambda$\th row-sum
  is $n!/\prod_{\text{chord $c$ of $\lambda$}} |c|$, as indicated in the following
  figure. 
  The $q$-analogue is $$\sum_\mu
  q^{|\lambda/\mu|/2}f_{\lambda/\mu}(q^{1/2}) =
  \frac{n!_q}{\prod_{\text{chord $c$ of $\lambda$}} |c|_q}.$$
  [A row specifies the lower path and allows any upper path above it.]
\begin{center}\includegraphics{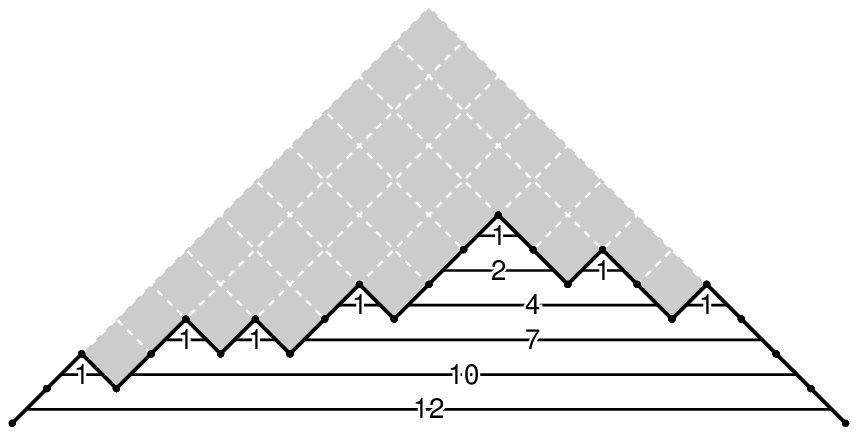}\end{center}
 $$\displaystyle \Rightarrow \frac{12!}{12\times 1\times 10\times 1\times 1\times 7\times1\times4\times1\times2\times1\times1} = 71280$$
\end{conjecture}

This formula holds whenever $n\leq 8$, and presumably in general.
An $n=4$ example is
\begin{center}
\psfrag{0}{$q^0$}
\psfrag{1}{$q^1$}
\psfrag{2}{$q^2$}
\psfrag{3}{$q^3$}
\psfrag{4}{$q^4$}
\psfrag{5}{$q^5$}
\psfrag{6}{$q^6$}
\includegraphics[scale=0.8]{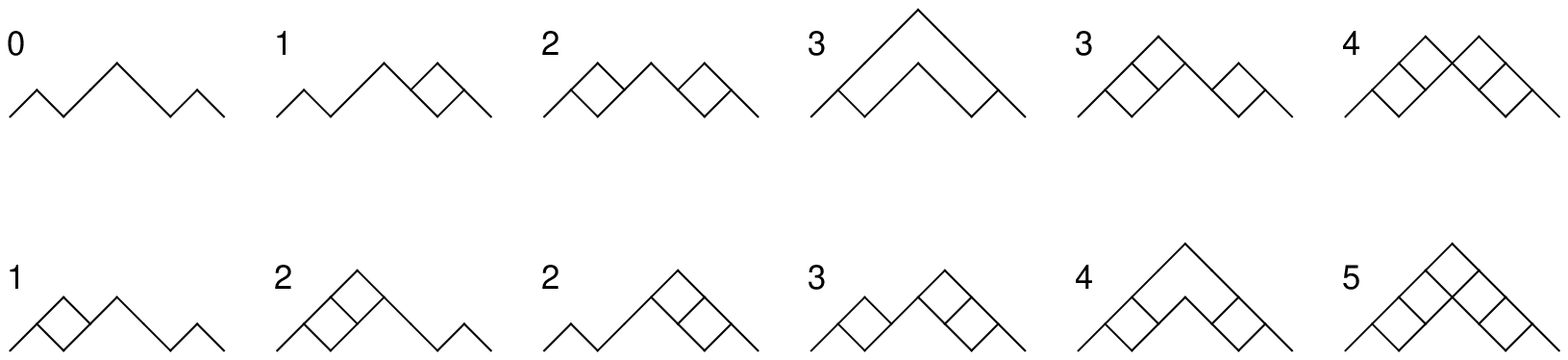}
\end{center}
$$\Rightarrow 1 + 2 q + 3 q^2 + 3 q^3 + 2 q^4 + q^5 = \frac{1\times (1+q)\times(1+q+q^2)\times (1+q+q^2+q^3)}{1\times (1+q)\times 1\times 1}$$

It is possible to prove this conjecture for the first row, corresponding
to the zigzag path, for which the row sum is $n!$.
This follows from the formula for the entries of the first row (\pref{row1}),
together with a very nice bijection between suitably labeled Dyck paths and
permutations due to de M\'edicis and Viennot \cite{MR1288802}.
The generating function for the $q$-analogue of the first row sum is the $q$-analogue
of a formula of Euler:
$$\sum_{n=0}^{\infty} n!\,x^n = \cfrac{1}{1-\cfrac{x}{1-\cfrac{x}{1-\cfrac{2x}{1-\cfrac{2x}{1-\cfrac{3x}{\ddots}}}}}}$$
the $q$-analogue of which is
$$\sum_{n=0}^{\infty}
 n!_q\,x^n = \cfrac{1}{1-\cfrac{1_qx}{1-\cfrac{q1_qx}{1-\cfrac{q2_qx}{1-\cfrac{q^22_qx}{1-\cfrac{q^23_qx}{\ddots}}}}}}$$

\old{
The recurrence is
$$
T_{x,y} = (q^{(y-1)/2} T_{x-1,y-1} + q^{y/2} T_{x-1,y+1})\times\begin{cases} 1 & \text{if $x$ and $y$ even} \\ (\frac{y+1}{2})_{q^{-1}} & \text{if $x$ and $y$ odd.} \end{cases}
$$
 The formula is
$$ T_{x,y} = q^{\lfloor y/2\rfloor^2/2} \Big\lceil \frac{x}{2} \Big\rceil !_q \binom{\big\lfloor \frac{x}{2} \big\rfloor}{\big\lfloor \frac{y}{2} \big\rfloor}_q.$$

Recall that $\binom{a}{b}_q = \binom{a-1}{b-1}_q + q^b\binom{a-1}{b}_q$.

If $x$ and $y$ are even we have
\begin{align*}
 T_{x,y} &= q^{(y-1)/2}T_{x-1,y-1} + q^{y/2}T_{x-1,y+1} \\ &= 
 q^{y/2-1/2}q^{y^2/8-y/2+1/2}\frac{x}{2} !_q \binom{\frac{x}{2} -1}{\frac{y}{2}-1}_q
+
 q^{y/2}q^{y^2/8}\frac{x}{2} !_q \binom{\frac{x}{2} -1}{\frac{y}{2}}_q\\
&=
 q^{y^2/8} \frac{x}{2} !_q \binom{\frac{x}{2}}{\frac{y}{2}}_q.
\end{align*}
If $x$ and $y$ are odd we have
\begin{align*}
 T_{x,y} &= \Big(\frac{y+1}{2}\Big)_{q^{-1}}
 (q^{(y-1)/2}T_{x-1,y-1} + q^{y/2}T_{x-1,y+1}) \\&= 
\Big(\frac{y+1}{2}\Big)_q
 \left[q^{(y-1)^2/8}           \frac{x-1}{2} !_q \binom{\frac{x-1}{2}}{\frac{y-1}{2}}_q +
       q^{1/2} q^{(y+1)^2/8}   \frac{x-1}{2} !_q \binom{\frac{x-1}{2}}{\frac{y+1}{2}}_q \right]\\
&=
q^{(y-1)^2/8} \Big(\frac{y+1}{2}\Big)_q \frac{x-1}{2} !_q \binom{\frac{x+1}{2}}{\frac{y+1}{2}}_q\\
&=
q^{(y-1)^2/8} \frac{x+1}{2} !_q \binom{\frac{x-1}{2}}{\frac{y-1}{2}}_q.
\end{align*}

}

\begin{conjecture} \label{columnsum}
For any column~$\mu$ of $M^{-1}$, the absolute values of the entries add up
to the product of the heights of the chords under $\mu$, as indicated in the figure.
The $q$-analogue is $$\sum_\lambda f_{\lambda/\mu}(q) = \prod_{\text{chord $c$ of $\lambda$}} (\text{height of $c$})_q.$$
[A column specifies the upper path and allows any lower path below it and above the zigzag.]
\begin{center}\includegraphics{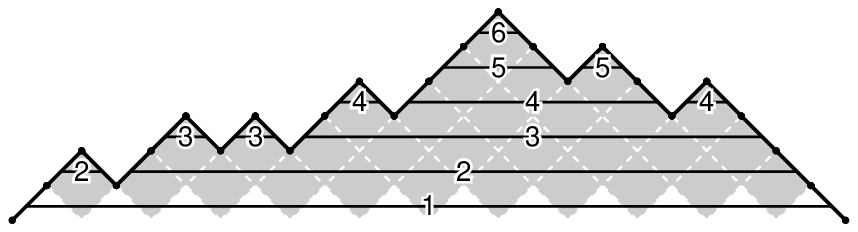}\end{center}
$$\Rightarrow 1\times2\times2\times3\times3\times3\times4\times4\times4\times5\times5\times6=1036800$$
\end{conjecture}

This formula holds whenever $n\leq 8$, and presumably in general.
An $n=4$ example is
\begin{center}
\psfrag{0}{$q^0$}
\psfrag{1}{$q^1$}
\psfrag{2}{$q^2$}
\psfrag{3}{$q^3$}
\psfrag{4}{$q^4$}
\psfrag{5}{$q^5$}
\psfrag{6}{$q^6$}
\includegraphics[scale=0.8]{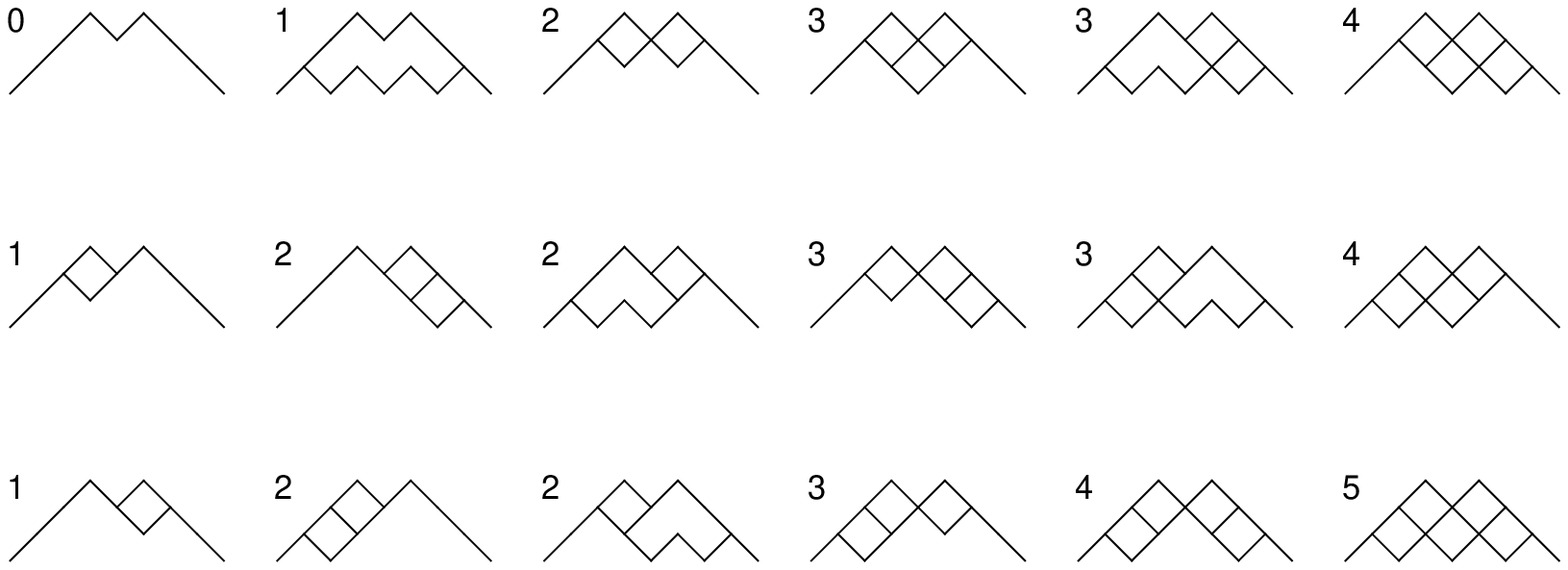}
\end{center}
$$\Rightarrow 1 + 3 q + 5 q^2 + 5 q^3 + 3 q^4 + q^5 = 1\times (1+q)\times(1+q+q^2)\times (1+q+q^2)$$

\section{Double-dimer marginals} \label{sec:ddimer}

We show how the previous results apply to the double-dimer
model.

Let $G$ be a finite bipartite planar graph with edges having
positive real weights, 
and $\No$ a 
set of $2n$ 
distinguished vertices on its outer face, which for simplicity
we assume alternate in color. 

\begin{window}[0,r,\psfrag{1}{$1$}\psfrag{2}{$2$}\psfrag{3}{$3$}\psfrag{4}{$4$}\includegraphics{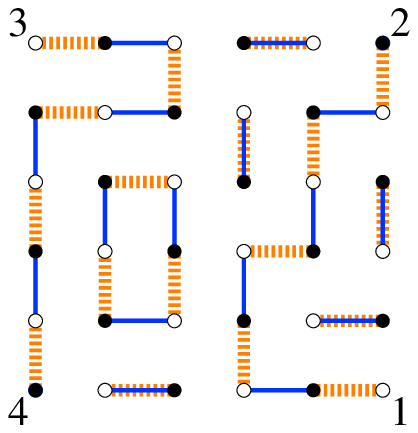},{}]
The {\bf double-dimer model with nodes $\No$} is the probability
measure on configurations obtained by superposing a random dimer cover
(perfect matching) of $\G$ with a random dimer cover of $\G\setminus\No$
(see figure at right).
Here by random dimer cover we mean a dimer cover chosen randomly for
the probability measure assigning each cover a probability
proportional to the product of edge weights of its constituent dimers.
Let $Z$ and $Z_{\No}$ denote the corresponding partition sums.
\end{window}
In such a configuration the nodes are joined in pairs with
lines of dimers. In \cite{KW:polynomial} we showed how to compute
the pairing probabilities, that is, the probability that a
given pairing of the nodes arises in a random configuration.
These pairing probabilities are rational functions of the
boundary measurements $\X_{i,j}$,
where $\X_{i,j}=Z_{i,j}/Z$ where $Z_{i,j}$ 
is the weighted sum of 
dimer covers of $\G\setminus\{i,j\}$  and $Z$ 
is the weighted sum
of dimer covers of $\G$.  See \cite{KW:polynomial}.
For example, the probability of the pairing $\{\{1,2\},\{3,4\}\}$,
which we write as $\,{}^1_2\!\mid\!{}^3_4$,
is $X_{1,2} X_{3,4} / (X_{1,2} X_{3,4} + X_{1,4} X_{2,3})$.

Knowing how to compute the full distribution on the
double-dimer pairings, we wish to compute its marginals. 
Suppose for example that we simply wish to know the
probability that node $i$ is paired with node $j$ when there is a very
large number of nodes, or infinitely many nodes.  This is an example
of a marginal probability question, and computing these marginal
probabilities can be non-trivial even when the full probability
distribution is known.  In this section we see how to compute these
marginal probabilities, and at the end we study a natural example that
has infinitely many nodes.

\subsection{Contiguous marginals}

We start by computing the probabilities of subpairings where the nodes
of the subpairing are contiguous.

We recall the definition of the quantities $D_S$ for a confining subset $S$ of nodes from
\cite[\S~3]{KW:polynomial}:
$$D_S=\det\left[(1_{i,j\in S}+1_{i,j\notin S})\times(-1)^{(|i-j|-1)/2}
\X_{i,j}\right]^{i=1,3\dots,2n-1}_{j=2,4,\dots,2n},$$
and the fact that
$$ \frac{D_S}{D_\varnothing} = \sum_\pi M_{S,\pi} \Pr(\pi).$$

This allows us to determine all contiguous marginals.
For example, row $\,{}^1_2\!\mid\!{}^3_4\!\mid\!{}^5_6\,$ of $M^{-1}$ yields the marginal
$$ \Pr(\,{}^1_2\!\mid\!{}^3_4\!\mid\!{}^5_6\!\mid\!\cdots) = \frac{D_{\{1,2,3,4,5,6\}} -D_{\{1,2,3,6\}} -D_{\{1,4,5,6\}} + D_{\{1,6\}}  -2 D_{\{1,3,4,6\}}}{D_\varnothing}.$$

\begin{window}[0,r,%
\psfrag{+D(1,2,3,4,5,6)}{$+D_{\{1,2,3,4,5,6\}}$}%
\psfrag{-D(1,2,3,6)}{$-D_{\{1,2,3,6\}}$}%
\psfrag{-D(1,4,5,6)}{$-D_{\{1,4,5,6\}}$}%
\psfrag{+D(1,6)}{$+D_{\{1,6\}}$}%
\psfrag{-D(1,3,4,6)}{$-D_{\{1,3,4,6\}}$}%
\includegraphics{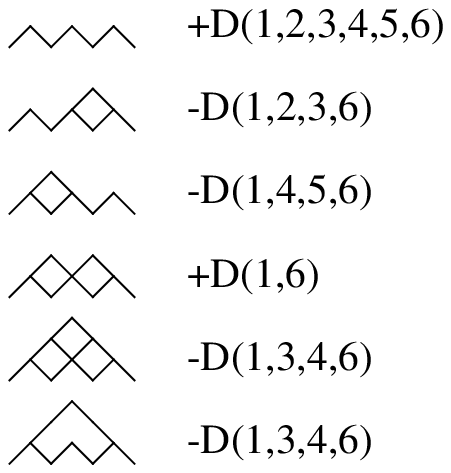},{}]
To derive this formula from the Dyck tilings, we would translate
$\,{}^1_2\!\mid\!{}^3_4\!\mid\!{}^5_6\,$ to a Dyck path \includegraphics{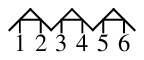}, which will be
the lower path of skew shapes for which we find cover-inclusive Dyck
tilings.  Each tiling corresponds to a term in the formula.  The sign
is negative when the area of the skew shape is odd.  The indices are
the locations of the odd-up and even-down steps of the upper boundary
of the shapes.
\end{window}

\subsection{Local non-contiguous marginals}

Next we consider marginal probabilities where the set of nodes
contained within the chords of the marginal are not contiguous, for
example
$\,{}^1_2\!\mid\!{}^3_6\!\mid\!{}^5_4\!\mid\!{}^{11}_{16}\!\mid\!\cdots$.
In this example there are three contiguous blocks of nodes contained
within chords: 1--6, 11, and 16.  We will show how to recursively
compute these marginals in terms of marginals with fewer contiguous blocks.
The previous subsection handled the base case of one contiguous block.

First observe that if a chord connects two non-contiguous blocks, then there is
a finite number of ways to ``fill in'' the space between the two connected blocks,
and each of these filled in marginals have fewer contiguous blocks.  In the
example above we have
$$
\Pr(\,{}^1_2\!\mid\!{}^3_6\!\mid\!{}^5_4\!\mid\!{}^{11}_{16}\!\mid\!\cdots) =
\Pr(\,{}^1_2\!\mid\!{}^3_6\!\mid\!{}^5_4\!\mid\!{}^{11}_{16}\!\mid\!{}^{13}_{12}\!\mid\!{}^{15}_{14}\!\mid\!\cdots) +
\Pr(\,{}^1_2\!\mid\!{}^3_6\!\mid\!{}^5_4\!\mid\!{}^{11}_{16}\!\mid\!{}^{13}_{14}\!\mid\!{}^{15}_{12}\!\mid\!\cdots)
$$

Next suppose that each chord connects only nodes within its block.
Using the method of the previous section, for each contiguous block we
associate a formal linear combination of subsets of nodes from that
block.  For example, in the first term above
$\,{}^1_2\!\mid\!{}^3_6\!\mid\!{}^5_4\!\mid\!{}^{11}_{16}\!\mid\!{}^{13}_{12}\!\mid\!{}^{15}_{14}\!\mid\!\cdots$,
we have $\{1,2,3,6\}-\{1,6\}+\{1,3,4,6\}$ for the first block and
$\{11,16\}-\{11,13,14,16\}$ for the second block.  We may extend the union operator $\cup$ linearly,
and take the union over all blocks of these formal linear combinations of subsets of nodes.  In this case we get
\begin{multline*}
 (\{1,2,3,6\}-\{1,6\}+\{1,3,4,6\}) \bigcup (\{11,16\}-\{11,13,14,16\}) = \\
   \{1,2,3,6,11,16\}-\{1,6,11,16\}+\{1,3,4,6,11,16\} \\
  -\{1,2,3,6,11,13,14,16\}+\{1,6,11,13,14,16\}-\{1,3,4,6,11,13,14,16\}.
\end{multline*}
Let $\sum_S \alpha_S S$ denote this resulting formal linear combination.
For any planar pairing $\pi$, we have
$$\sum_S \alpha_S M_{S,\pi} = \begin{cases}1 & \text{$\pi$ is consistent with the marginal of interest}\\0 & \text{$\pi$ connects each node of a block to another node of the same block,}\\[-2pt] &\text{but not consistent with marginal}\\ 0 & \text{$\pi$ connects a node of some block to a node not in a block }\\ \text{integer} & \text{otherwise.}\end{cases}
$$
But pairings $\pi$ that fall into the fourth case above can be
``filled in'' as discussed above, and the resulting marginals have
fewer blocks and can thus recursively be expressed in terms of the
$D_S$'s.  These can then be subtracted from $\sum_S \alpha_S S$.

\subsection{Nonlocal marginals}

\begin{figure}[t]
\centerline{\includegraphics[width=.5\textwidth]{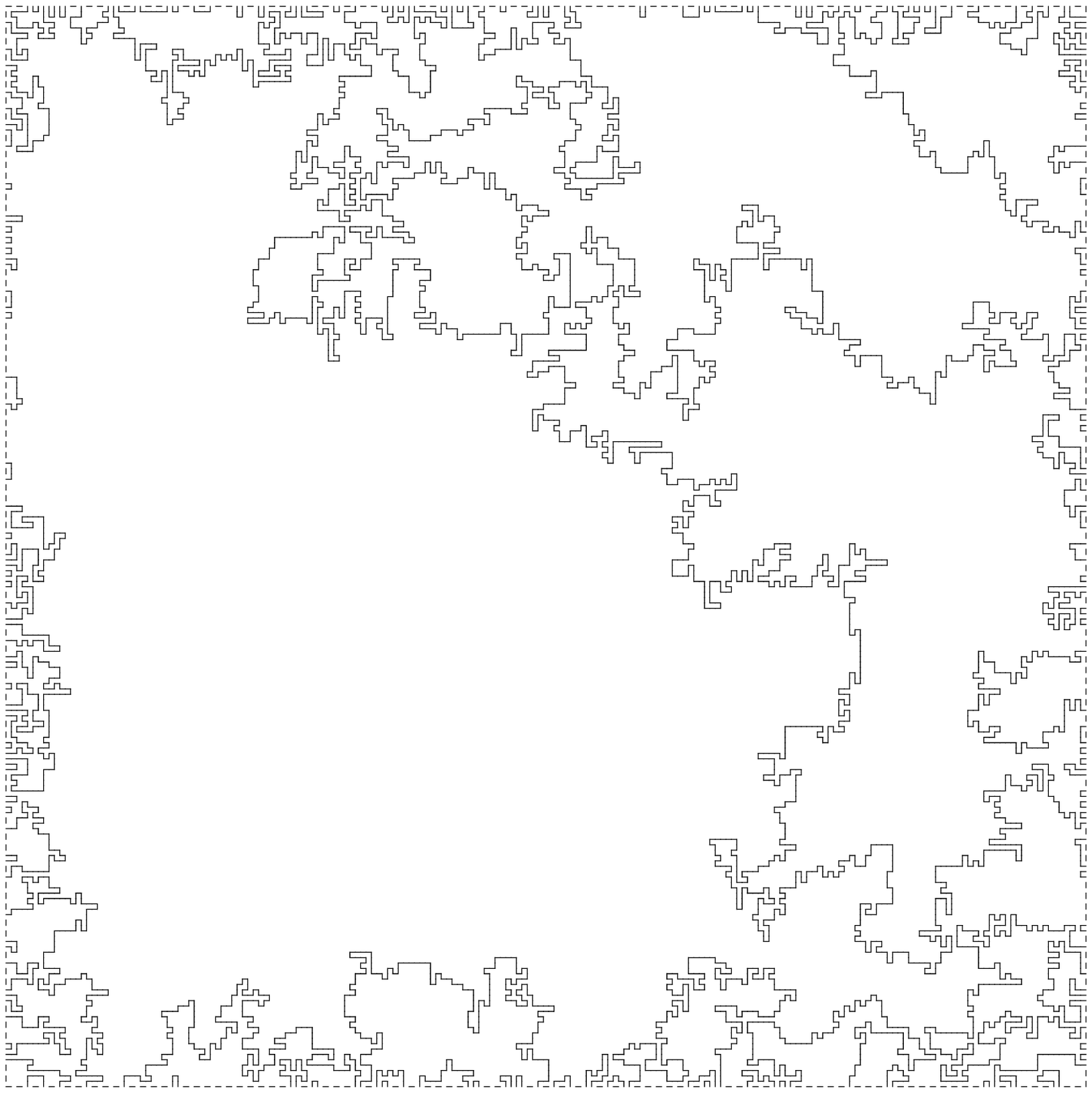}\hfil
\includegraphics[width=.45\textwidth]{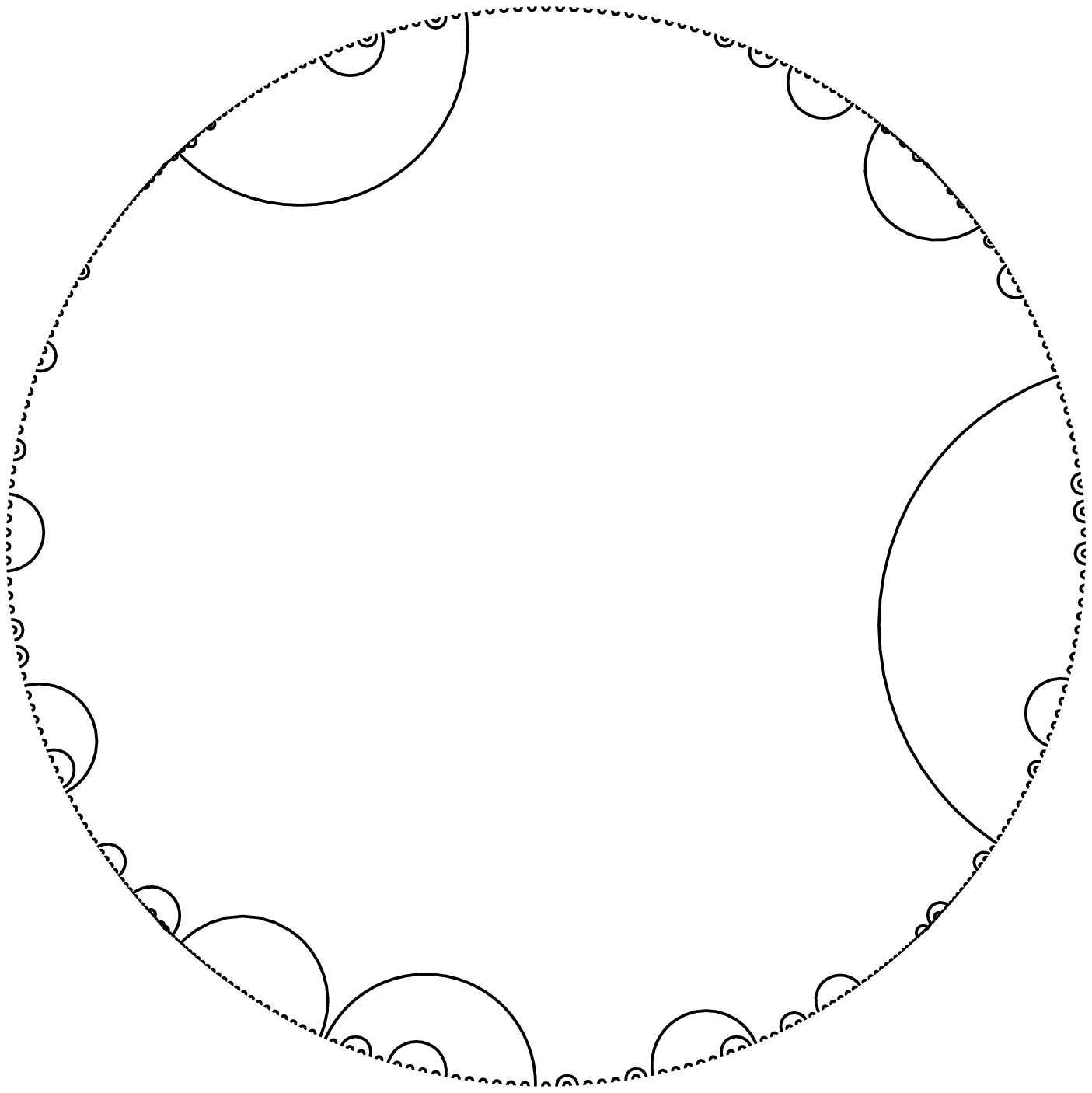}}
\caption{On the left are the chains in a random double-dimer configuration in a box-shaped region with many nodes all along its boundary.
  The loops and doubled edges of the double-dimer configuration are not shown here.
  In the figure on the right, the box was conformally mapped to the disk, and the chains connecting
  points on the circle are replaced with arcs.  We expect that in the scaling limit,
  such random chord diagrams will be M\"obius invariant.
  \label{boxpair}}
%\centerline{\includegraphics[width=.5\textwidth]{arcs}}
%  \label{arcs}}
\end{figure}

We have seen how to compute local marginal pairing probabilities.  If
we wish to understand the scaling limit of double-dimer pairings such
as the ones in \fref{boxpair}, we need to be able to compute
marginal pairing probabilities where the nodes in question are well
separated.  For example, we'd like to be able to compute the
probability that the upper left corner is paired with the middle of
the left edge while the upper right corner and lower right corner are
paired with each other.  This is a marginal probability on just four
nodes, but there are many nodes separating these four nodes.  
Unfortunately we do not know how to compute these
non-local marginals in general, even for four nodes out of $2n$.
Experiment seems to indicate that the $4$-node
marginals depend on all the $X_{i,j}$.

\subsection{Evenly spaced nodes}

\begin{figure}[htbp]
\includegraphics[width=\textwidth]{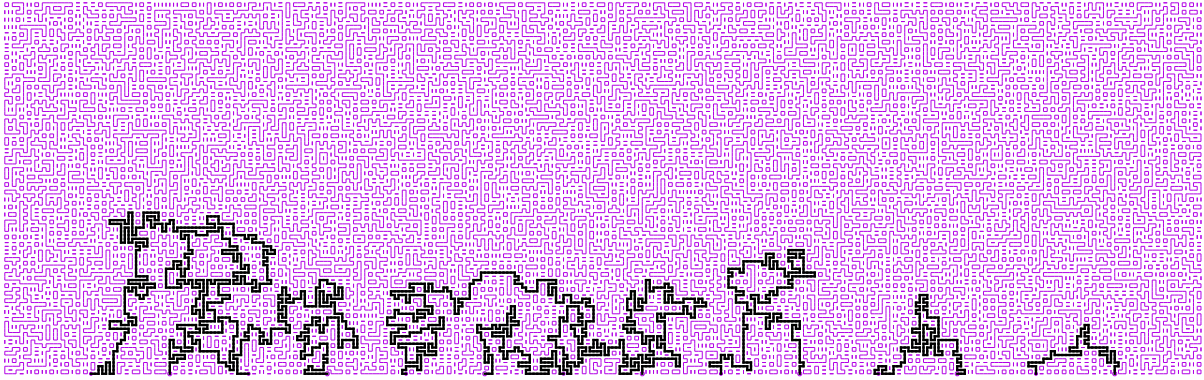}
\caption{A random double-dimer configuration with nodes of
  alternating color at evenly spaced locations.  We are interested in
  the limit where the domain converges to the upper half plane, the lattice
  spacing tends to $0$, and there is a node at each integer.
  \label{ddimer-zchain}}
\end{figure}

As an application of these ideas, consider the double-dimer model on
$\eps\Z\times\eps\N$, and its scaling limit on $\R\times\R^+$, and
suppose that there are $2n$ boundary nodes which alternate in color and are
at locations $x_1,\dots,x_{2n}$ (see \fref{ddimer-zchain}).
In the scaling limit $X_{i,j}=1/|x_i-x_j|$, (see \cite{MR1782431}),
and we have the remarkable formula
  $$\frac{D_S}{D_\varnothing} = \prod_{i\in S, j\notin S} |x_j-x_i|^{(-1)^{1+i+j}}$$
\cite[Lemma~5.2]{KW:polynomial}.
In \cite{KW:polynomial} we proved this formula for finite numbers of
nodes on the real line bounding the upper half plane, but it also
holds for a finite number of nodes on a circle bounding a disk.  We
can use the right-hand side to define $D_S/D_\phi$ for non-balanced
sets $S$ of the nodes.  Using
$$
\frac{D_{\{i\}}}{D_\varnothing} = \prod_{j\neq i} |x_j-x_i|^{(-1)^{1+i+j}},
$$
we can rewrite
$$
\frac{D_S}{D_\varnothing} = \prod_{i\in S} \frac{D_{\{i\}}}{D_\varnothing} \prod_{\substack{i,j\in S\\i<j}} |x_i-x_j|^{2(-1)^{i+j}}.
$$

Now suppose that the domain is the unit disk and there are $2n$ nodes,
one at each of the $2n$\th roots of unity.  Let $\zeta=e^{2\pi i/n}$ and $\omega=e^{\pi i/n}$.
$$
\prod_{j=1}^n (z-\zeta^j) = z^n -1 = (z-1) (z^{n-1}+z^{n-2}+\dots+z+1)
$$
Thus
 $$
\frac{D_{\{i\}}}{D_\varnothing} = \frac{\prod_{j=1}^n |\omega-\zeta^j|}{\prod_{j=1}^{n-1} |1-\zeta^j|} = \frac{2}{n}.
$$
We also have
$$
|x_i-x_j| = 2\left|\sin\frac{\pi(i-j)}{2n}\right|.
$$
Hence
$$
\frac{D_S}{D_\varnothing} = \frac{1}{n^{|S|}} \frac{\displaystyle\prod_{\substack{i,j\in S\\i<j\\\text{$j-i$ even}}} \sin^2\frac{\pi(i-j)}{2n}}{\displaystyle\prod_{\substack{i,j\in S\\i<j\\\text{$j-i$ odd}}} \sin^2\frac{\pi(i-j)}{2n}}.
$$
For example, if $j-i$ is odd then
$$ \frac{D_{\{i,j\}}}{D_\varnothing} = \frac{4}{n^2} \frac{1}{|x_i-x_j|^2} = \frac{1}{n^2} \frac{1}{\sin^2\frac{\pi(i-j)}{2n}} \to \frac{4/\pi^2}{(i-j)^2}$$
in the $n\to\infty$ limit.  More generally, for any finite balanced set $S$, in the $n\to\infty$ limit
$$
\frac{D_S}{D_\varnothing} = \left(\frac{2}{\pi}\right)^{|S|} \frac{\displaystyle\prod_{\substack{i,j\in S\\i<j\\\text{$j-i$ even}}} (i-j)^2}{\displaystyle\prod_{\substack{i,j\in S\\i<j\\\text{$j-i$ odd}}} (i-j)^2}.
$$

For example, the probability that $1$ pairs with $2$ and $3$ pairs with $4$ is
 $$\Pr[\cdots\!\mid\!{}^1_2\!\mid\!{}^3_4\!\mid\!\cdots] = \frac{D_{\{1,2,3,4\}}}{D_\varnothing} - \frac{D_{\{1,4\}}}{D_\varnothing} = \left(\frac{2}{\pi}\right)^4 \frac{16}{9} - \left(\frac{2}{\pi}\right)^2 \frac{1}{9} = 0.246979\dots.$$

\section{Pairings in groves}

\begin{window}[0,r,\psfrag{1}{$1$}\psfrag{2}{$2$}\psfrag{3}{$3$}\psfrag{4}{$4$}\includegraphics{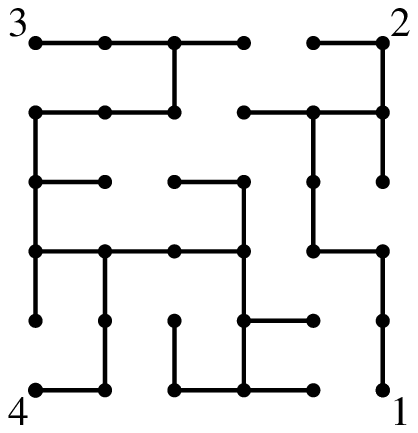},{}]
Many results for dimers have analogues for trees and vice versa (see
e.g., \cite{MR633646,KPW,MR2099145}).  The analogue of the
double-dimer pairings of the previous section are ``grove pairings.''
Given a graph (not necessarily bipartite) with a specified subset $\No$ of the vertices called nodes,
a grove is a forest such that each tree contains at least one of the
nodes (see figure at right).
A grove pairing is a grove in
which each tree contains exactly two such nodes.  If all such groves
have probability proportional to the product of their edge weights, we
are interested in the probability that a random grove has a specified
pairing type.
\end{window}

\hspace{\parindent}
An important formula for grove pairings, which is due to Curtis,
Ingerman, and Morrow \cite{MR1657214} (see also Fomin
\cite{MR1837248}), expresses grove pairing probabilities in terms of
determinants of a submatrix of the ``(electrical) response matrix''.
If node $i$ is held at $1$ volt while the other nodes are held at $0$
volts, then $L_{i,j}$ denotes the current flowing out of the network
through node $j$.  Though not obvious from this definition,
$L_{i,j}=L_{j,i}$ (see e.g., \cite{MR1652692}).  
If there are $2n$ nodes $\{v_1,\dots,v_n,w_1,\dots,w_n\}$,
the Curtis-Ingerman-Morrow formula is equivalent to
$$
\det\begin{bmatrix}L_{i,j}\end{bmatrix}^{i=v_1,\ldots,v_n}_{j=w_1,\ldots,w_n} = \sum_{\sigma\in S_n} (-1)^\sigma \frac{Z_{v_1,w_{\sigma(1)}|\cdots|v_n,w_{\sigma(n)}}}{Z_{v_1|\cdots|v_n|w_1|\cdots|w_n}},
$$
where $Z_\pi$ is the weighted sum of groves where the nodes are
connected according to the partition $\pi$.  This formula holds for
any graph, whether or not it is planar.  In the case of planar graphs,
we require the nodes to be a subset of the vertices on the outer face,
which we number in circular order.  In this case, if all the $v_i$'s
are contiguous, then there is only one pairing between the $v$'s and
$w$'s for which the term on the right-hand-side is nonzero, so for
this pairing there is a simple determinant formula.

When the graph is planar and the nodes are on the outer face, there
are $C_n = \binom{2n}{n}\frac{1}{n+1}$ noncrossing pairings $\pi$
of the nodes, and $\binom{2n}{n}/2$ formulas relating the $Z_\pi$'s
coming from the CIM formula, so Dub\'edat conjectured that all of the
$Z_\pi/Z_{1|\cdots|2n}$'s are determined by the CIM determinants
\cite{MR2253875}.  Here we show how to adapt our double-dimer results
from the previous section to prove this conjecture.  Essentially the
same formulas hold as for the double-dimer pairings, except that there
are some extra minus signs in the formulas for grove pairings.

Let $S^*=\{v_1,\ldots,v_n\}\subset\{1,\dots,2n\}$.  Then the CIM
determinant enumerates (with sign) grove pairings compatible with
$S^*$, where a pairing is compatible with $S^*$ iff every chord
connects an element of $S^*$ to an element not in $S^*$.  Given $S^*$,
we can define $S$ by
$$S = \text{(odds of $S^*$)} \cup \text{(evens not in $S^*$)}.$$
Then the pairing is compatible with $S$ iff no chord connects an
element of $S$ to an element not in $S$, which is precisely the same
notion of compatibility that we used in the previous section on
double-dimer pairings.

We can rewrite the previous equation as
$$
\det\begin{bmatrix}L_{i,j}\end{bmatrix}^{i\in S^*}_{j\notin S^*} = \sum_{\pi} (-1)^{\sigma(\pi,S^*)} M_{S,\pi} \frac{Z_\pi}{Z_{1|\cdots|2n}},
$$
since whenever $Z_\pi\neq 0$, we have $M_{S,\pi}=1$.  Here
$\sigma(\pi,S^*)$ is the permutation mapping $S^*$ to the
complement of $S^*$, when these sets are listed in sorted order.
Conveniently, this sign $(-1)^{\sigma(\pi,S^*)}$ does not depend on $S^*$,
so we can write it as $(-1)^\sigma$:

\begin{lemma}
  Let $\pi$ be a planar pairing.  Suppose $S^*$ is a set of nodes that
  is compatible with $\pi$ in the sense that each chord connects a
  node of $S^*$ to a node of $(S^*)^c$.  If the nodes of $S^*$ are
  arranged in sorted order and likewise for $(S^*)^c$, then the
  permutation defined by the pairing has sign which is determined by
  $\pi$ alone and is independent of $S^*$.
\end{lemma}
\begin{proof}
  It suffices to compare the signs for $S_1^*$ and $S_2^*$ which
  differ at one chord $\{a,b\}$ of the pairing, where say $a<b$ and
  $a\in S_1^*$ and $b\in S_2^*$.  When $a\in S_1^*$ is replaced with
  $b$, it is moved to the right a number of times for the list to
  remain in sorted order, and this number equals the number of chords
  $\{c,d\}$ where $a<c<d<b$.  Likewise when $b\in S_1^*$ is replaced
  with $a$, it is moved to the left the same number of times for the
  list to remain in sorted order.  Thus the parity of the permutation
  connecting $S_1^*$ to $(S_1^*)^c$ changes an even number of times when
  it is transformed to the permutation connecting $S_2^*$ to $(S_2^*)^c$.
\end{proof}

Therefore we can solve for the $Z_\pi$'s in the same way that we
solved for the $\Pr[\pi]$'s in the previous section:
\begin{theorem}
If $\pi$ is a noncrossing pairing of $\{1,\dots,2n\}$, for groves in a planar graph we have
$$
(-1)^\pi \frac{Z_\pi}{Z_{1|\cdots|2n}} = \sum_{S} M^{-1}_{\pi,S} \,\det \begin{bmatrix}L_{i,j}\end{bmatrix}^{i\in S^*}_{j\notin S^*}.
$$
\end{theorem}

For example, when there are $6$ nodes, the partition function for the
pairing $\,{}^1_2\!\mid\!{}^3_4\!\mid\!{}^5_6\,$ is given by
\begin{multline*}
\frac{Z_{{}^1_2\mid{}^3_4\mid{}^5_6}}{Z_{1|2|3|4|5|6}} =
\begin{vmatrix}
 L_{1,2} & L_{1,4} & L_{1,6} \\
 L_{3,2} & L_{3,4} & L_{3,6} \\
 L_{5,2} & L_{5,4} & L_{5,6}
\end{vmatrix}
-\begin{vmatrix}
 L_{1,2} & L_{1,5} & L_{1,6} \\
 L_{3,2} & L_{3,5} & L_{3,6} \\
 L_{4,2} & L_{4,5} & L_{4,6} 
\end{vmatrix}
\\
-\begin{vmatrix}
 L_{1,3} & L_{1,4} & L_{1,6} \\
 L_{2,3} & L_{2,4} & L_{2,6} \\
 L_{5,3} & L_{5,4} & L_{5,6} 
\end{vmatrix}
+\begin{vmatrix}
 L_{1,3} & L_{1,5} & L_{1,6} \\
 L_{2,3} & L_{2,5} & L_{2,6} \\
 L_{4,3} & L_{4,5} & L_{4,6} 
\end{vmatrix}
-2\begin{vmatrix}
 L_{1,4} & L_{1,5} & L_{1,6} \\
 L_{2,4} & L_{2,5} & L_{2,6} \\
 L_{3,4} & L_{3,5} & L_{3,6} 
\end{vmatrix}
.
\end{multline*}

\begin{window}[0,r,%
\psfrag{+D(1,2,3,4,5,6)}{$+\det\begin{bmatrix}L_{i,j}\end{bmatrix}^{i=1,3,5}_{j=2,4,6}$}%
\psfrag{-D(1,2,3,6)}{$-\det\begin{bmatrix}L_{i,j}\end{bmatrix}^{i=1,3,4}_{j=2,5,6}$}%
\psfrag{-D(1,4,5,6)}{$-\det\begin{bmatrix}L_{i,j}\end{bmatrix}^{i=1,2,5}_{j=3,4,6}$}%
\psfrag{+D(1,6)}{$+\det\begin{bmatrix}L_{i,j}\end{bmatrix}^{i=1,2,4}_{j=3,5,6}$}%
\psfrag{-D(1,3,4,6)}{$-\det\begin{bmatrix}L_{i,j}\end{bmatrix}^{i=1,2,3}_{j=4,5,6}$}%
\includegraphics{tile-row-dets}\ \ \,\,,{}]
To derive this formula from the Dyck tilings, as with the double-dimer pairings,
we translate
$\,{}^1_2\!\mid\!{}^3_4\!\mid\!{}^5_6\,$ to a Dyck path \raisebox{-2pt}{\includegraphics{pairing}},
 which will be
the lower path of skew shapes for which we find cover-inclusive Dyck
tilings (see figure to right).  Each tiling corresponds to a term in the formula.  The sign
is negative when the area of the skew shape is odd.  The rows are indexed
by the locations of the up steps of the upper boundary, and the columns
are indexed by the down steps of the upper boundary.  In addition there is
a sign on the left-hand-side which is the sign of the permutation associated
with the pairing.
\end{window}

\begin{window}[0,r,%
\psfrag{+D(1,6)}{$+\det\begin{bmatrix}L_{i,j}\end{bmatrix}^{i=2,3,5}_{j=4,6,1}$}%
\psfrag{-D(1,3,4,6)}{$-\det\begin{bmatrix}L_{i,j}\end{bmatrix}^{i=2,3,4}_{j=5,6,1}$}%
\includegraphics{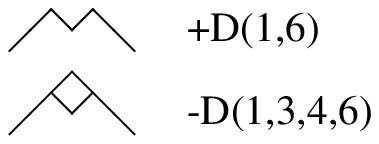}\ \ \,\,,{}] If we cyclically
rotate the indices, then the pairing
$\,{}^1_2\!\mid\!{}^3_4\!\mid\!{}^5_6\,$ translates to the Dyck
path \raisebox{-2pt}{\includegraphics{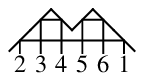}}, which allows us to write
an equivalent formula with just two determinants, as shown in the figure to the right.
\end{window}

Note that we only know how to compute the grove partition function
$Z_\pi$ for a complete pairing $\pi$; unlike the situation for double-dimer pairings,
we do not know how to compute marginal
probabilities.

\phantomsection
\pdfbookmark[1]{References}{bib}
\bibliographystyle{hmralpha}
\bibliography{bc}

\end{document}